\documentclass[10pt]{article}

\usepackage{auto-pst-pdf}
\usepackage{amsmath, latexsym, amsfonts, amssymb, amsthm}
\usepackage{graphicx, color,hyperref,dsfont,epsfig,tikz,caption,wrapfig,subfig,manfnt}
 \usepackage{relsize}
   \newcommand{\myWarning}[1]{\textcolor{red}{\textbf{\relscale{0.5}\textdbend\relscale{2}\,#1}}}
\usepackage{pstricks}
\usepackage{datetime}
\usepackage{movie15}
\hypersetup{
    linktoc=page,
    linkcolor=red,          
    citecolor=blue,        
    filecolor=blue,      
    urlcolor=cyan,
    colorlinks=true           
}

\numberwithin{equation}{section}

\newtheorem{theorem}{Theorem}[section]
\newtheorem{lemma}[theorem]{Lemma}
\newtheorem{proposition}[theorem]{Proposition}

\theoremstyle{definition}

\renewcommand{\tilde}{\widetilde}          
\DeclareMathSymbol{\leqslant}{\mathalpha}{AMSa}{"36} 
\DeclareMathSymbol{\geqslant}{\mathalpha}{AMSa}{"3E} 
\DeclareMathSymbol{\eset}{\mathalpha}{AMSb}{"3F}     
\renewcommand{\leq}{\;\leqslant\;}                   
\renewcommand{\geq}{\;\geqslant\;}                   

\newcommand{\C}{\mathbb{C}}
\newcommand{\R}{\mathbb{R}}

\newcommand{\E}{\mathds{E}}
\renewcommand{\P}{\mathds{P}}

\newcommand{\hf}{\frac{_1}{^2}}
\renewcommand{\log}{\ln}

\usepackage{xparse}

\DeclareDocumentCommand \Pmp { m m o} {
\IfNoValueTF{#3}
{P_{#1}^{#2}}
{P_{#1}^{#2}\left(#3\right)}
}

\DeclareDocumentCommand \Emp { m m o} {
\IfNoValueTF{#3}
{E_{#1}^{#2}}
{E_{#1}^{#2}\left[#3\right]}
}

\DeclareDocumentCommand \Pbr { m m m m o } {
\IfNoValueTF{#5}
{P_{#1}^{#2\stackrel{#4}{\rightarrow} #3}}
{P_{#1}^{#2\stackrel{#4}{\rightarrow} #3}\left(#5\right)}
}

\DeclareDocumentCommand \Ebr { m m m m o } {
\IfNoValueTF{#5}
{E_{#1}^{#2\stackrel{#4}{\rightarrow} #3}}
{E_{#1}^{#2\stackrel{#4}{\rightarrow} #3}\left[#5\right]}
}

\def\S{\mathbb{S}}

\def\bi{\begin{itemize}}
\def\ei{\end{itemize}}

\def\bnum{\begin{enumerate}}
\def\enum{\end{enumerate}}
\def\<#1{\langle #1 \rangle}

\title{Renormalizability of Liouville Quantum Gravity  at the Seiberg bound}

\author{ Fran\c{c}ois David \footnote{Institut de Physique Th\'eorique,
CNRS, URA 2306, CEA, IPhT, Gif-sur-Yvette, France.}, Antti Kupiainen \footnote{University of Helsinki, Department of Mathematics and Statistics, P.O.
Finland. Supported by the Academy of Finland}, R\'emi Rhodes \footnote{Universit{\'e} Paris-Est Marne la Vall\'ee, LAMA, Champs sur Marne, France.} \footnotetext[2]{Partially supported by grant ANR-11-JCJC  CHAMU.},
 Vincent Vargas \footnote{ENS Ulm, DMA, 45 rue d'Ulm,  75005 Paris, France.} }

\begin{document}

\maketitle

\begin{abstract}
Liouville Quantum Field Theory can be seen as a probabilistic theory of 2d Riemannian metrics $e^{\phi(z)}dz^2$, conjecturally describing scaling limits of discrete $2d$-random surfaces. The law of the random field $\phi$ in LQFT depends on weights $\alpha\in \R$ that  in classical Riemannian geometry parametrize power law singularities in the metric. 
A rigorous construction of LQFT    has been carried out in \cite{DKRV} in the case when the weights  are below the so called Seiberg bound: $\alpha<Q$ where $Q$  parametrizes the random surface model in question. These correspond to conical singularities in the  classical setup. In this paper, we construct  LQFT  in the case when  the Seiberg bound is saturated which can be seen as the probabilistic version  of  Riemann surfaces with cusp singularities.  Their construction involves methods from Gaussian Multiplicative Chaos theory at criticality. 
\end{abstract}

\begin{center}
\end{center}
\footnotesize



\noindent{\bf Key words or phrases:}  Liouville Quantum Gravity, quantum field theory, Gaussian multiplicative chaos, KPZ formula, KPZ scaling laws, Polyakov formula, punctures, cusp singularity, uniformization theorem.


\normalsize

\tableofcontents


\section{Introduction and main results}


 In this paper we continue the rigorous study of the two dimensional  Liouville Quantum Theory  or Liouville Quantum Gravity (LQG) started in \cite{DKRV} by means of probabilistic tools. We summarize in this introduction only the main points established  there and then we describe our results. We  refer the reader  to  \cite{cf:Da,DistKa,cf:KPZ,nakayama,Pol,seiberg}  for seminal physics references  on LQG and to \cite{DKRV} for  more background and references in maths and physics.
 
 \subsection{Summary of LQG}
 
 The LQG is a probabilistic theory of Riemannian metrics of the form $e^{\gamma X} g$
 where $ g$ is a fixed smooth Riemannian metric on a two dimensional surface $\Sigma$ and $X$ is a random field on $\Sigma$  whose law is formally given in terms of a  functional integral 
\begin{equation}\label{pathintegral}
\E[F(X)]=Z^{-1}\int F(X)e^{-S_L(X,{g})}DX
\end{equation}
where $Z$ is a normalization constant, $DX$ stands for a formal uniform measure on  some space of maps  $X:\Sigma\to\R$ and 
\begin{equation}\label{actionfalse}
S_L(X,{g}):= \frac{1}{4\pi}\int_{\Sigma}\big(|\partial^{{g}}X |^2+QR_{{g}} X  +4\pi \mu e^{\gamma X  }\big)\,d\lambda_{{g}}  
\end{equation}
is the so-called Liouville action, where $\partial^{{g}}$, $R_{{g}} $ and $\lambda_{{g}}$ respectively stand for the gradient, Ricci scalar curvature and volume measure in the metric $ {g}$ (see Section \ref{sec:backgr} for the basic definitions used in here). The parameter $\mu> 0$ is called ``cosmological constant'',   $\gamma\in [0,2)$ is a parameter that is determined by the random surface model in question and  $Q=\frac{2}{\gamma}+\frac{\gamma}{2}$.

The rigorous definition of the integral \eqref{pathintegral} in the case of the sphere $\Sigma=\S^2$ was carried out in \cite{DKRV}. In the sequel, by stereographic projection, we will often identify $\S^2$ with $\R^2 \cup \lbrace \infty \rbrace$ (the identification should be clear from the context).  We take $g=g(z)|dz|^2$ conformally equivalent (Section \ref{sec:backgr}) to the standard round metric 
and set
\begin{equation}\nonumber
X=c+X_{{g}}
\end{equation}
where  $X_{ g}$ is the Gaussian Free Field  with vanishing $\lambda_{ g}$-mean on the sphere $\S^2$, $\int X_g\lambda_{{g}} =0$ (Section \ref{sec:backgr}), and $c\in\R$ . Definition of the Liouville term $e^{\gamma X  }d\lambda_{{g}} $ 
requires a regularization and a renormalization. Let  
\begin{equation}\label{law}
X_{g,\epsilon}(x)=\frac{1}{2\pi}\int_0^{2\pi}X_g(x+\epsilon e^{i\theta})\,d\theta.
 \end{equation}
be the circle average regularization of the GFF and define the random measure
\begin{equation}\label{circleav}
dM_{\gamma}=\lim_{\epsilon\to 0}\epsilon^{\frac{\gamma^2}{2}}e^{\gamma (X_{{{g}},\epsilon}+Q/2 \ln {g})}  \,d\lambda,
 \end{equation}
where $d \lambda$ is the standard Lebesgue measure. The limit is in probability in the sense of weak convergence of measures. 
$M_{\gamma}$ is a random measure on $\S^2$ with total mass $M_{\gamma}(\S^2)$ almost surely finite.

The precise definition of $e^{-S_L}DX$ is given by the measure
\begin{equation}\label{lioumeas}
d\mu_L(X):=e^{- \frac{_Q}{^{4\pi}}\int_{\S^2}R_{{g}} X d\lambda_{{g}} - \mu e^{\gamma c}M_{\gamma}(\S^2)}d\P(X_{ g})\,dc
 \end{equation}
defined on $H^{-1}(\S^2, \lambda_{{g}})\times\R$. This measure is not normalizable to a probability measure due to the $c$ variable. Indeed, by the Gauss-Bonnet theorem, $ \frac{_Q}{^{4\pi}}\int_{\S^2}R_{{g}}  d\lambda_{{g}}=2Q$ and thus the total mass of
$d\mu_L$ is infinite due to divergence of  the $c$-integral at $-\infty$.
This divergence is related in an interesting way to the (quantum) geometry described by the
LQG. 

First, Liouville theory is {\it conformally invariant}. With this respect, one should distinguish the classical and quantum theory. 
On the classical level, one studies the functional   \eqref{actionfalse} on a functional space and chooses $Q=\frac{2}{\gamma}$ to ensure that it is invariant under the following action of group of
M\"obius transformations $\psi$ of $\S^2$
\begin{equation}\label{liouinv}
S_L(X_\psi,g)=S_L(X,g) 
 \end{equation}
 where $X_\psi=X\circ\psi+\frac{_1}{^\gamma}\log \frac{g_\psi}{g}$, $g_\psi=|\psi'|^2g\circ\psi$ \footnote{ In fact, following standard conventions in the physics literature, one can absorb the metric dependence of the action by shifting the field $X$, i.e. introducing the Liouville field $\phi(x)=X(x)+\frac{Q}{2} \log g(x)$. This convention was used in \cite{DKRV} and will also be used later in this paper: see definition \eqref{defLiouvillefield}. In this case, the action can be written with respect to the standard Euclidean distance as background metric on $\R^2$ and takes on the following form:
\begin{equation*}
S_L(\phi, g):= \frac{1}{4 \pi} \int_{\R^2} ( |\partial \phi|+4 \pi \mu e^{\gamma \phi}) d \lambda, \quad \phi \sim \frac{Q}{2} \ln g,
\end{equation*} 
where recall that $\lambda$ is the standard Lebesgue measure and $\sim$ says that the difference is bounded. In this case, we get the invariance $S_L(\phi \circ \psi + Q \log |\psi'|,g)=S_L(\phi,g) $ which is more familiar with the physics literature. One also gets an analogue on the quantum side.
}. On the quantum level, one is interested in the path integral formulation \eqref{pathintegral}. It was proven in   \cite{DKRV} that
 \begin{equation}\label{lioumeas}
\int F(X_\psi)d\mu_L(X)=\int F(X)d\mu_L(X)
 \end{equation}
for $F\in L^1(\mu_L)$ where now  $X_\psi=X\circ\psi+\frac{_Q}{^2}\log \frac{g_\psi}{g} $
  and $Q=\frac{2}{\gamma}+\frac{\gamma}{2}$ \footnote{The extra term $\frac{\gamma}{2}$ comes from the $\epsilon^{\frac{\gamma^2}{2}}$-renormalization of \eqref{circleav}.}.

The invariance of $\mu_L$ under the action of the non-compact group $SL(2,\C)$ is one reason to expect that $\mu_L$ has infinite mass. The other reason has to do with the fact that the classical action functional \eqref{actionfalse} is not bounded from below. Indeed the minimizers of the Liouville action are  the solutions of the classical  {\it Liouville equation} $R_{e^{\gamma X}g}=-2\pi\mu\gamma^2$. These are metrics conformally equivalent to $g$  with  constant negative curvature. On $\S^2$ there are no smooth negative curvature metrics since by the Gauss-Bonnet theorem the total curvature is positive. 

Both of these problems can be fixed by adding punctures to the sphere. Let $z_i\in\S^2$, $\alpha_i\in\R$, $i=1,\dots,n$ and consider the measure $\prod_ie^{ \alpha_iX(z_i)}\mu_L$. The  {\it vertex operators} $e^{\alpha X(z)}$ require renormalization. Define
\begin{equation}\label{vertex}
V_{\alpha,\epsilon}(z)=\epsilon^{\frac{\alpha^2}{2}}e^{\alpha (c+X_{{{g}},\epsilon}(z)+Q/2 \ln {g}(z))}
\end{equation}
and set
\begin{equation}\label{vertexcorr}
\Pi_{\bf \alpha, z,\epsilon}:= \int\prod_iV_{\alpha_i,\epsilon}(z_i)d\mu_L(X),
\end{equation}
 Due to the  $c$ variable this exists only if $\sum_i\alpha_i>2Q$. In \cite{DKRV} it was shown that $\Pi_{\bf \alpha, z}=\lim_{\epsilon\to 0}\Pi_{\bf \alpha, z,\epsilon}>0$  if and
only if  $\alpha_i<Q$ for all $i$. Briefly, the reason for this is as follows. One can absorb the vertex operators in \eqref{vertexcorr} by an application of the Girsanov transform i.e. by a shift of the gaussian field $X_{g}\to X_{ g}+H$ with
\begin{equation}\label{Hdef}
H(z)=\sum_i\alpha_iG_{g}(z,z_i).
\end{equation}
where $G_{  g}$ is the covariance of the GFF. Taking $g$ equal to the round metric $\hat g$ (with $R_{\hat g}=2$), this leads to 
\begin{equation}\label{vertexcorr1}
\Pi_{\bf \alpha, z}=K({\bf z}) \int_\R e^{(\sum_i\alpha_i-2Q) c }  (\E  e^{ - \mu e^{\gamma c} \int e^{\gamma H}dM_\gamma })dc
\end{equation}
with
\begin{equation}\label{Kdef}
K({\bf z})=\prod_i\hat{g}(z_i)^{-\frac{\alpha_i^2}{4}+\frac{Q}{2}\alpha_i}e^{\hf \sum_{i\not =j}\alpha_i\alpha_{j}G_{\hat{g}}(z_i,z_{j})+\frac{\ln 2-1/2}{2} \sum\alpha_i^2}
\end{equation}
Due to the logarithmic singularity of $G_{\hat g}$ the integrand $e^{\gamma H(z)}$ blows up as $ |z-z_i|^{-\alpha_i\gamma}$ when $z\to z_i$.  
By analyzing the modulus of continuity  of the Gaussian multiplicative chaos measure  it was shown in \cite{DKRV} that $\int e^{\gamma H}dM_\gamma $
is a.s. finite if and only if $\alpha_i<Q$. Moreover,   the probability measures
\begin{equation}\label{vertexcorr3}
\P_{\bf \alpha, z,\epsilon}:=\Pi_{\bf \alpha, z,\epsilon}^{-1} \prod_iV_{\alpha_i,\epsilon}(z_i)\, d \mu_L(X).
\end{equation}
converge
\begin{equation}\label{vertexcorr4}
\lim_{\epsilon\to 0}\P_{\bf \alpha, z,\epsilon}=\P_{\bf \alpha, z}.
\end{equation}

The bounds  $\sum_i\alpha_i>2Q$ and   $\alpha_i<Q$ are called the {\it Seiberg bounds} (see also \cite{seiberg}). They lead to the conclusion that to have a nontrivial correlation function  of vertex operators one needs at least {\it three} of them. Note that fixing three points on the sphere fixes also  the $SL(2,\C)$ invariance. When we have three insertions with $\alpha_i=\gamma$ the law of the chaos measure $e^{\gamma c}M_\gamma$ under $\P_{\bf \alpha, z}$ is conjectured to agree with the scaling limit of planar maps decorated with a critical statistical mechanics model and conformally embedded onto the sphere. In  \cite{DKRV} this correspondence was checked explicitly for the law of the total volume $e^{\gamma c} M_\gamma(\S^2)$.

It is instructive to consider the classical problem in the presence of vertex operators i.e. with their logarithms subtracted from the Liouville action (the reader may consult \cite{troyanov} for further explanations than those given below). Performing the substitution
$
X=c+H+\tilde X
$ with  $\int \tilde Xd\lambda_{\hat g}=0$ this functional becomes  using the classical value $Q=\frac{2}{\gamma}$
\begin{equation}\label{punctures1}
(\frac{_4}{^\gamma}-\sum_i\alpha_i)c+ \frac{_1}{^{4\pi}}\int_{\S^2}|\partial^{}\tilde X |^2+\mu e^{\gamma c}\int_{\S^2} e^{\gamma \tilde X  }\,d\lambda_{{ g}}  -\log K({\bf z})
\end{equation}
where $g=e^{\gamma H}\hat g$ and $K({\bf z})$ is given in \eqref{Kdef}. \eqref{punctures1} is  bounded from below iff $\sum_i\alpha_i>\frac{4}{\gamma}$. 
The volume form $\lambda_g$ is integrable provided $\alpha_i<\frac{2}{\gamma}$. 
These are the classical versions of the Seiberg bounds. The metric $g$ has negative curvature in the complement of the punctures and  a conical singularity at $z_i$ with the angle $\theta_i$ of the cone equal to $\pi(2-\alpha_i\gamma)$, that is, a neighborhood of $z_i$ is asymptotically isometric to a neighborhood of the tip of a cone with the map $z\to z^{\frac{\theta}{2\pi}}$. The second Seiberg bound thus states that the angle needs to be positive. In this case \eqref{punctures1} has a bounded minimizer $(c,\tilde X)$. Note also that again we need at least three punctures to satisfy the Seiberg bounds.

\subsection{Main Results}

 In this paper we extend the analysis of   \cite{DKRV} to the case of vertex operators with weight $Q$, called $Q$-punctures. Classically  as $\alpha_i\to\frac{2}{\gamma}$ the minimizer $\tilde H$ of \eqref{punctures1} is no more bounded and the metric $e^{\gamma\tilde H(z)}g$  blows up  as (take $z_i=0$)  $(|z|\log|z|)^{-2}|dz|^2$. Geometrically this corresponds to a cusp singularity of the metric with the volume around the puncture finite whereas distances become infinite. These are the so called parabolic solutions of  the  Liouville equation. 
 
 In the quantum case   as noted above, $\Pi_{\bf \alpha, z}$ tends to zero as $\alpha_i\to Q$. However,
 a simple renormalization suffices for obtaining a nontrivial limit
  \footnote{\myWarning{}  In fact we will use a  slightly different definition of $\Pi_{\bf \alpha, z,\epsilon}$ from \eqref{vertexcorr} because we will regularize simultaneously the vertex operators \eqref{vertex} and the measure $\mu_L$ defined by \eqref{lioumeas} (see equation \eqref{eq:defPiGB}).}:
 
\begin{theorem}\label{th:seiberg1} 
Let $\sum_i\alpha_i>2Q$ and $\alpha_i\leq Q$ with exactly $k$ of the $\alpha_i$ equal to $Q$. Then the limit
\begin{equation}\label{limit1}
\lim_{\epsilon\to 0}({-\log\epsilon})^{\frac{k}{2}} \Pi_{\bf \alpha, z,\epsilon}:= \Pi_{\bf \alpha, z}
\end{equation}
exists and is strictly positive. Moreover, the limit
 \begin{equation}\label{limit}
\lim_{\epsilon\to 0}\P_{\bf \alpha, z,\epsilon}:=\P_{\bf \alpha, z}
\end{equation} 
exists in the sense of weak convergence of measures on $H^{-1}(\S^2)$.
\end{theorem}

This theorem means that the vertex operator $V_Q$ needs an additional factor $({-\log\epsilon})^\hf $ for its normalization in addition to the Wick ordering used for $\alpha<Q$.  
This normalization is familiar from the Seneta-Heyde normalization needed for the construction of the multiplicative chaos measure at criticality \cite{Rnew7,Rnew12}. As in that context an important ingredient in the proof of convergence \eqref{limit1} is to show that the limit agrees (up to a multiplicative constant) with the one constructed with the {\it derivative vertex operator}
 \begin{equation}\label{dervertexQ}
\tilde V_{Q,\epsilon}(z)=-\frac{d}{d\alpha}|_{\alpha=Q}V_{\alpha,\epsilon}(z)=-(Q\log\epsilon+c+X_{g,\epsilon}+\frac{_Q}{^2}\log g)V_{\alpha,\epsilon}(z).
\end{equation} 
Let $ \tilde \Pi_{\bf \alpha, z,\epsilon}$ be the correlation function where for $\alpha_i=Q$ we use $\tilde V_{Q,\epsilon}(z)$. Then

 \begin{theorem}\label{th:seiberg2}  
\begin{equation}\label{limit2}
\lim_{\epsilon\to 0} \tilde \Pi_{\bf \alpha, z,\epsilon}= (\frac{_\pi}{^2})^\frac{k}{2}\Pi_{\bf \alpha, z}
\end{equation}
\end{theorem}

The convergence \eqref{limit} extends to  functions of the chaos measure. Let $\E_{\bf \alpha, z,\epsilon}$ denote expectation in $\Pi_{\bf \alpha, z,\epsilon}$ and
let $F=F(X,\nu)$ be a bounded continuous function on $H^{-1}(\S^2)\times {\mathcal M}(\S^2)$ where ${\mathcal M}(\S^2)$ denotes Borel measures on $\S^2$. Define the {\it Liouville measure}
\begin{equation}\label{defLiouvillemeas}
Z:=e^{\gamma c}M_\gamma
\end{equation}
and the {\it Liouville field}
\begin{equation}\label{defLiouvillefield}
\phi:=X+\frac{_Q}{^2}\log g.
\end{equation}
Then
 
\begin{theorem}\label{th:seiberg3} With the assumptions of Theorem \ref{th:seiberg1},
$\E_{\bf \alpha, z,\epsilon}F(\phi, Z)$ converges as $\epsilon\to 0$ to a  limit  
$\E_{\bf \alpha, z}F$ which  is conformally covariant
$$
\E_{\bf \alpha, z}F(\phi, Z)=\E_{\bf \alpha, \psi(z)}F(\phi\circ\psi+ Q \ln |\psi'|, Z \circ \psi)
$$
and independent of $g$ in the conformal equivalence class of $\hat g$. 
Moreover, the law of $Z(\S^2)$ under
$\P_{\bf \alpha, z}$  is given by  
the Gamma distribution 
\begin{equation}\label{gammadist}
E_{\bf \alpha, z}F(Z(\S^2))=\frac{\mu^{\frac{\sigma}{\gamma}}}{\Gamma({\frac{\sigma}{\gamma}})}\int_0^{\infty}F(y)y^{{\frac{\sigma}{\gamma}-1}}e^{-\mu y}\,dy,\ \ \ \sigma:=\sum_i \alpha_i- 2Q
\end{equation}
and the law of the random measure $Z(\cdot)/A$ conditioned on $Z(\S^2)=A$ does  not depend on $A$.
\end{theorem}

\noindent {\bf Remark 1.} {\it The  correlation functions $\Pi_{\bf \alpha, z}$ have the same properties as in the $\alpha_i<Q$ case proven in  \cite{DKRV}: {\bf conformal covariance, Weyl covariance and KPZ scaling}. Since the statements are identical we refer the reader to \cite{DKRV} recalling here only the KPZ formula
for the $\mu$-dependence:}
$$
\Pi_{\bf \alpha, z}=
\mu^{\frac{2Q-\sum_i\alpha_i}{\gamma}}\Pi_{\bf \alpha, z}|_{\mu=1} .
$$
\vskip 2mm

\noindent {\bf Remark 2.} {\it Theorem \ref{th:seiberg3}  is the quantum analog of the convergence of the elliptic solutions of the Liouville equation to parabolic ones as one saturates the second Seiberg bound. As in the classical case, the quantum volume of the metric is a.s. finite.
}
\vskip 2mm

\noindent {\bf Remark 3.} {\it
With some extra work it should be possible to prove that the measures $\P_{\bf \alpha, z}$ with $\alpha_i<Q$ for all $i=1,\dots n$ converge as $\alpha_i\uparrow Q$, $i=1,\dots k$ to the $\P_{\bf \alpha, z}$ constructed in this paper by proving that
  \begin{equation}\label{limit2}
\lim_{\alpha_i\uparrow Q}\prod_{i=1}^k(Q-\alpha_i)^{-1}\Pi_{\bf \alpha, z}
\end{equation}
has a limit.  We leave that question as an open problem.  }

 \vskip 2mm

\noindent {\bf Remark 4.} {\it  It is natural to ask about the convergence of the quantum laws $\P_{\alpha,z}$ to the classical solutions of the Liouville equation i.e. the {\bf semiclassical limit} $\gamma\to 0$. For this, let us take, for $i=1,\dots,k$ $\alpha_i=Q$ and for $i>k$
$$\alpha_i=\frac{\chi_i}{\gamma}$$
with $\chi_i<2$
and $\mu=\frac{\mu_0}{\gamma^2}$  for some constant $\mu_0>0$. Then we conjecture that the law of $\gamma X$ under $\P_{\alpha,z}$ converges towards the minimizer $c+\tilde X$ of equation
\eqref{punctures1} which has  cusp singularities at  $z_i$, $i\leq k$ and conical ones at the remaining $z_i$.}

 
 \vskip 3mm
 Finally, let us mention that these $Q$-punctures are especially important to understand how to embed conformally onto the sphere random planar maps with spherical topology weighted by a $c=1$ conformal field theory (like the Gaussian Free Field). Indeed, in the case $c=1$, one can   formulate the conjecture developed in \cite[subsection 5.3]{DKRV} with $\gamma=2$ and $Q=2$: the vertex operators with $\gamma=2$ in \cite[conjecture 2]{DKRV} are precisely the $Q$-punctures constructed in this paper.  Though we do not treat explicitly  the case $\gamma=2$, the techniques we develop adapt to this case.


\section{Background and notations}\label{sec:backgr}
Here we recall some background, taken from \cite{DKRV}, which will be used throughout the paper.

\medskip \noindent {\bf Basic notations.} $B(x,r)$ stands for the ball centered at $x$ with radius $r$. We let $\bar C({\R^2})$ stand for the space of continuous functions on $\R^2$ admitting a finite limit at infinity. In the same way, $\bar C^k({\R^2})$ for $k\geq 1$ stands for the space of $k$-times differentiable functions  on $\R^2$ such that all the derivatives up to order $k$ belong to $\bar C({\R^2})$. 

\medskip\noindent  {\bf Metrics on  $\S^2$.}  
The  sphere $\S^2$ can be mapped by stereographic projection to the plane which we view both as   ${\R^2}$ and as $\C$. 
We take as the  background metric  the  round metric on $\S^2$ which becomes on $\R^2$ and on $\C$
$$\hat{g}=\frac{4}{(1+|x|^2)^2}dx^2=\frac{2}{(1+\bar zz)^2}(dz\otimes d\bar z+d\bar z\otimes dz).$$
Its Ricci scalar curvature is $R_{\hat g}=2$. 
 We say a metric $g=g(x)dx^2$ is conformally equivalent to $\hat g$ if 
$$g(x)=e^{\varphi(x)}\hat{g}(x)
$$ 
with $\varphi\in \bar C^2({\R^2})$  such that $\int_{\R^2}|\partial \varphi|^2\,d\lambda<\infty$. We often identify the metrics $g$ with their densities $g(x)$ (or $g(z)$) with respect to the Euclidean metric. We denote the volume measure $g(x) \lambda(dx)$ by $ \lambda_g(dx)$ where $\lambda$ is the Lebesgue measure on $\R^2$. The total volume of $\hat g$ is $\int_{\R^2}d\lambda_{\hat g}=4\pi$.

Given any metric $g$ conformally equivalent to the spherical metric, we let
$H^1(\S^2)=H^1(\R^2,g)$ be the Sobolev space defined as the closure of $\bar C^\infty({\R^2})$ with respect to the Hilbert-norm
\begin{equation}\label{Hnorm}
\int_{\R^2}|h|^2\,d\lambda_{g}+\int_{\R^2}|\partial h|^2\,d\lambda.
\end{equation}

\medskip \noindent {\bf Gaussian free fields. }
For each metric $g$ conformally equivalent to $\hat{g}$, we consider a Gaussian Free Field $X_g$ with vanishing $\lambda_g$-mean on the sphere,  that is a centered  Gaussian random distribution with covariance kernel given by the Green function $G_g$ of the problem
\begin{equation}\nonumber
\triangle_g u=-2\pi f \quad\text{on }\R^2,\quad  \int_{\R^2} u \,d\lambda_g =0
\end{equation}
i.e. $u=\int G_g(\cdot,z)f(z)\lambda_g(dz):=G_gf$. 
In case of  the round metric, we have the explicit formula 
\begin{equation}\label{hatGformula}
\E X_{\hat g}(z)X_{\hat g}(z')= G_{\hat g}(z,z')=\ln\frac{1}{|z-z'|}-\frac{1}{4}(\ln\hat g(z)+\ln\hat g(z'))+\ln 2-\hf.
\end{equation}
The GFF $X_g$ lives almost surely   in the dual space $H^{-1}(\S^2)$ of $H^{1}(\S^2)$,  and this space does not depend on the choice of the metric $g$ in the conformal equivalence class of $\hat{g}$. 

%
%

   
\medskip \noindent {\bf Gaussian multiplicative chaos. }
The circle average regularization \eqref{circleav} $X_{\hat g,\epsilon}$ of the free field $X_{\hat g}$ satisfies\begin{equation}\label{circlegreen}
\lim_{\epsilon\to 0}\E[X_{\hat{g},\epsilon} (x)^2]+\ln \epsilon
+\frac{1}{2}\ln \hat{g}(x)=  \ln 2-\hf
\end{equation} uniformly on $\R^2$. Define now the measure
\begin{equation}\label{Meps}
M_{\gamma,\epsilon}:=\epsilon^{\frac{\gamma^2}{2}}e^{\gamma (X_{{\hat{g}},\epsilon}+Q/2 \ln \hat{g})}  \,d\lambda.
\end{equation}
For $\gamma\in [0,2)$,   we have the convergence  in probability
\begin{equation}\label{law}
M_{\gamma}=\lim_{\epsilon\to 0}M_{\gamma,\epsilon}
=e^{\frac{\gamma^2}{2}(\ln 2-\hf)} \lim_{\epsilon\to 0}e^{\gamma X_{\hat{g},\epsilon}-\frac{\gamma^2}{2} \E[ X_{\hat{g},\epsilon}^2] } \,d\lambda_{\hat{g}}
 \end{equation}
in the sense of weak convergence of measures. This limiting measure is non trivial and is a (up to a multiplicative constant) Gaussian multiplicative chaos \cite{cf:Kah,review} of the field $X_{\hat{g}}$ with respect to the measure $\lambda_{\hat{g}}$.

\section{Partition of the probability space}\label{sec:main}

The singularity at the $Q$-punctures will be studied by partitioning the probability space according to the maximum of the circle average fields around them.  As we will see this is a local operation and  it will suffice to consider the case with only one $Q$-insertion, say $\alpha_1=Q$, $\alpha_i<Q$, $i>1$. Also, 
 we will work from now on with the round metric $\hat g$; the general case $g=e^{\phi}\hat g$ is treated as in \cite{DKRV}. It will be convenient to modify the definition  \eqref{vertexcorr} slightly around the $Q$-insertion.  For this remove an  $\epsilon $-disc around the $z_1$
$$D_\epsilon:=\R^2\setminus B(z_1,\epsilon) $$ 
and define

\begin{equation}\label{eq:defPiGB}
 \Pi_{\bf \alpha, z,\epsilon}(F)=  \int_\R e^{\sigma c}  \E\Big[F( c+X_{\hat{g}}) \prod_i 
 V_{z_i,\alpha_i,\epsilon} (z_i) e^{-\mu
 e^{\gamma c}M_\gamma(D _\epsilon)}  \Big]\,dc
\end{equation}
where we use throughout the paper the notation 
\begin{equation}\label{sigmadef}
\sigma:=\sum_i\alpha_i-2Q
\end{equation}
 as in \eqref{gammadist}. We have then 
\begin{equation}\label{eq:defPiGB1}
 \E_{z,\alpha,\epsilon}F=  \Pi_{\bf \alpha, z,\epsilon}(F)/ \Pi_{\bf \alpha, z,\epsilon}(1).
 \end{equation}
 The Girsanov argument then gives $ \Pi_{\bf \alpha, z,\epsilon}(F)= K_\epsilon({\bf z})A_{\epsilon}(F)$ with
 \begin{equation}\label{eq:defPiGB2}
 A_{\epsilon}(F)=\int_\R e^{\sigma c}  \E\Big[F( c+X_{\hat{g}}+H_\epsilon) e^{-\mu
 e^{\gamma c}\int_{D _\epsilon}e^{\gamma H_\epsilon}dM_\gamma}  \Big]\,dc.
\end{equation}
where (recall  \eqref{Hdef})
\begin{equation}\label{def:Hepsilon}
H_\epsilon(z)= \sum_i\alpha_i\int_0^{2\pi}G_{\hat{g}}(z_i+\epsilon e^{i\theta},z)\frac{_{d\theta}}{^{2\pi}}
\end{equation}
and $K_\epsilon({\bf z})$ converges to $K$ of \eqref{Kdef} as $\epsilon\to 0$.

Similarly for the derivative vertex operator \eqref{dervertexQ}
we get
\begin{equation}\label{eq:defPiGB2}
 \tilde\Pi_{\bf \alpha, z,\epsilon}(F)= -K_\epsilon({\bf z}) \int_\R e^{\sigma c}  \E\Big[F( c+X_{\hat{g}}+H_\epsilon) (Q\log\epsilon+\tilde H_\epsilon+c+X_{\hat g,\epsilon}(z_1)+\frac{_Q}{^2}\log \hat{g}(z_1))e^{-\mu
 e^{\gamma c}\int_{D _\epsilon}e^{\gamma H_\epsilon}dM_\gamma}  \Big]\,dc.
\end{equation}
where
\begin{equation}\label{def:tildeHepsilon}
\tilde H_\epsilon=  \sum_i\alpha_i\int_0^{2\pi}\int_0^{2\pi}G_{\hat{g}}(z_1+\epsilon e^{i\theta_2},z_i+\epsilon e^{i\theta_2})\frac{_{d\theta_1}}{^{2\pi}}\frac{_{d\theta_2}}{^{2\pi}}.
\end{equation}
Using \eqref{hatGformula} we see that the $Q\log\epsilon$ singularity in \eqref{eq:defPiGB2} is cancelled by the one in the $i=1$ term in \eqref{def:tildeHepsilon} so that 
$
Q\log\epsilon+\tilde H_\epsilon+\frac{_Q}{^2}\log g(z_1)
$ is bounded uniformly in $\epsilon$. Since $\Pi_{z,\alpha,\epsilon}(F)\to 0$ as $\epsilon\to 0$ (\cite{DKRV}) we conclude that the limit, if it exits, of  $ \tilde\Pi_{z,\alpha,\epsilon}(F)$ equals the limit of $K_\epsilon({\bf z})\tilde A_{\epsilon}(F)$ where

\begin{equation}\label{eq:defPiGB22}
 \tilde A_{\epsilon}(F)= 
  \int_\R e^{\sigma c}  \E\Big[F( c+X_{\hat{g},\epsilon}+H_\epsilon) (-c-X_{\hat g})e^{-\mu
 e^{\gamma c}\int_{D _\epsilon}e^{\gamma H_\epsilon}dM_\gamma}  \Big]\,dc.
\end{equation}
 
Hence Theorems \ref{th:seiberg1} and  \ref{th:seiberg2} follow if we prove 

\begin{proposition}\label{2limits}  Let $F$ be bounded and continuous on $H^{-1}(\S^2)$. Then
the following limits 
\begin{equation}\label{two limits}A(F)=\lim_{\epsilon\to 0} (-\log\epsilon)^\hf A_{\epsilon}(F)=\frac{_2}{^\pi}\lim_{\epsilon\to 0} \tilde A_{\epsilon}(F)
\end{equation}
exist and $A(1)>0$. 
\end{proposition}

Now we partition the probability space according to the values of the maximum of the mapping $u\mapsto X_{\hat{g},u}(z_1)$ over $u\in [\epsilon,1]$.  So we set 
\begin{align}
M_{n,\epsilon}= &\big\{\max_{u\in[\epsilon,1]} X_{\hat{g},u}(z_1)   \in [n-1,n]\big\},\,\,\,\,n\geq 1,\label{def:mn}\\
M_{0,\epsilon}= &\big\{\max_{u\in[\epsilon,1]} X_{\hat{g},u}(z_1)  \leq 0\big\}.\label{def:m0},
\end{align}
and we expand the integral $ A_\epsilon(F)$ along the partition made up of these sets $(M_{n,\epsilon})_n$: 
\begin{equation}
 A_\epsilon(F) =\sum_{n\geq 0} \int_\R e^{\sigma c }  \,\E\Big[\mathbf{1}_{M_{n,\epsilon}}   F( c  +X_{\hat{g}} +H_{\epsilon})  e^{- \mu e^{\gamma c} \int_{D_\epsilon}e^{\gamma H_{\epsilon}} \,dM_\gamma }\Big]\,dc
 :=\sum_{n\geq 0} A_\epsilon(F,n).
\end{equation}
For $ \tilde A_\epsilon(F)$ we write
\begin{equation}\nonumber
 \tilde A_\epsilon(F) =\sum_{n\geq 0} (\tilde A_\epsilon(F,n)+ B_\epsilon(F,n))
\end{equation}
with
\begin{equation}\label{tildeAep}
 \tilde A_\epsilon(F,n)=\int_\R e^{\sigma c }  \,\E\Big[\mathbf{1}_{M_{n,\epsilon}} \big(  n-X_{\hat{g},\epsilon}(z_1)\big)  F( c  +X_{\hat{g}} +H_{\epsilon})  e^{- \mu e^{\gamma c} \int_{D_\epsilon}e^{\gamma H_{\epsilon}} \,dM_\gamma} \Big]\,dc
\end{equation}
and 
\begin{equation}\label{Bep}
 B_\epsilon(F,n)=-\int_\R e^{\sigma c }  \,\E\Big[\mathbf{1}_{M_{n,\epsilon}} (n+c) F( c  +X_{\hat{g}} +H_{\epsilon})  e^{- \mu e^{\gamma c} \int_{D_\epsilon}e^{\gamma H_{\epsilon} }\,dM_\gamma} \Big]\,dc.
\end{equation}
Note that $ \tilde A_\epsilon(F,n)\geq 0$ for $F\geq 0$.
We prove

\begin{lemma}\label{crucial1}  Let $F$ be bounded and continuous on $H^{-1}(\S^2)$. Then for all $n\geq 0$ the limits
\begin{equation}\label{two limits}
A(F,n)=\lim_{\epsilon\to 0} (-\log\epsilon)^\hf A_{\epsilon}(F,n)=\sqrt{2/\pi}\lim_{\epsilon\to 0} \tilde A_{\epsilon}(F,n).
\end{equation}
exist and $A(1,n)>0$. Moreover
\begin{align}\label{est2}
\sum_{n\geq 0}\sup_{\epsilon\in]0,1]} (-\log\epsilon)^\hf A_{\epsilon}(1,n)<\infty  \\
\label{est3}
\sum_{n\geq 0}\sup_{\epsilon\in]0,1]}\tilde A_{\epsilon}(1,n)<\infty \\
\sum_{n\geq 0}B_\epsilon(F,n)\to 0,\quad \text{ as }\epsilon\to 0.\label{est4}
\end{align}
\end{lemma}
Proposition \ref{2limits} then follows from Lemma \ref{crucial1} since $\lim_{\epsilon\to 0} A_{\epsilon}(F,\epsilon)=\sum A(F,n)$ follows from  \eqref{two limits} and \eqref{est2} by the dominated convergence theorem, idem for $\tilde A$. The remaining part of this paper is devoted to proving this lemma.

\section{Decomposition of the GFF and Chaos measure}
With no loss, we may set the $Q$ insertion at $z_1=0$ and suppose that the other $z_i$ are in the complement  of the disc $B(0,1)$. 
We further denote  by $\mathcal{F}_\delta$ ($\delta>0$) the sigma-algebra generated by the field $X_{\hat{g}}$ "away from the disc $B(0,\delta)$", namely
\begin{equation}
\mathcal{F}_\delta=\sigma\{X_{\hat{g}}(f);{\rm supp}\, f\in B(0,\delta)^c\}.
\end{equation}
$\mathcal{F}_\infty$ stands for the sigma algebra generated by $\bigcup_{\delta>0}\mathcal{F}_\delta$. 
 First observe that (see \cite{DKRV,DMS,review})
\begin{lemma}\label{BM}
For all $\delta>0$, the process 
$$ t\mapsto  X_{\hat{g},\delta e^{-t}}(0)-X_{\hat{g},\delta}(0)  $$
evolves as a Brownian motion independent of the sigma algebra $\mathcal{F}_\delta$.
\end{lemma}
The following decomposition of the field $X_{\hat{g}}$ will be crucial for the analysis:

\begin{lemma}\label{independence} The field $X_{\hat{g}}$  may be decomposed 
\begin{equation}
 X_{\hat{g}}(z)= X_{\hat{g},|z|}(0) + Y(z)
\end{equation} 
where 
the process $ r\in \R_+^*\mapsto X_{\hat{g},r}(0)$  is independent of the field $Y(z)$. The latter has the following covariance
\begin{equation*}
\E[    Y(re^{i\theta}) Y(r'e^{i\theta'}) ]  = \ln \frac{ r \vee r'}{|re^{i\theta}-r'e^{i\theta'}|}
\end{equation*}
\end{lemma}

\noindent {\it Proof.} From \eqref{hatGformula} we get using rotational invariance
$\E [ X_{\hat{g}}(z)X_{\hat{g},|z'|}(0)]=
\E [X_{\hat{g},|z|}(0)X_{\hat{g},|z'|}(0)].
$ 
which in turn leads to independence:
$$
\E X_{\hat{g}}(z)X_{\hat{g}}(z')=\E X_{\hat{g},|z|}(0)X_{\hat{g},|z'|}(0)+
\E Y(z)Y(z').
$$
Furthermore  we calculate
$$
\E[    Y(re^{i\theta}) Y(r'e^{i\theta'}) ]  = \ln {|re^{i\theta}-r'e^{i\theta'}|^{-1}}- \frac{1}{4 \pi^2}\int_{0}^{2 \pi}\int_{0}^{2 \pi}  \ln |re^{i u}-r'e^{i v}|du dv.
$$
The claim follows from  $\int_0^{2\pi}  \ln |re^{i \theta}-r'e^{i u}| du  =  \ln ( r \vee r')$.
\qed

Now, we get the decomposition
$$
M_\gamma(dz)= c_\gamma  \: \hat g(z)  |z|^{\frac{\gamma^2}{2}}e^{\gamma X_{\hat{g},|z|}(0)}\,M_\gamma(dz,Y)
$$
where $M_\gamma(dz,Y)$ is the multiplicative chaos measure of the field $Y$ with respect to the Lebesgue measure $\lambda$ (i.e. $\E\, M_\gamma(dz,Y)=\lambda(dz)$) and $c_\gamma:= e^{\frac{\gamma^2}{2} (\ln 2-1/2) -\frac{\gamma^2}{2}  \E[  X_{\hat{g},1}(0)^2  ] }$ is some constant.

We will now make change of variables $z=e^{-s+i\theta}$, $s\in \R_+$, $\theta\in [0,2\pi)$ and let $\mu_{Y}(ds,d\theta)$ be  the multiplicative chaos measure of the field $Y(e^{-s+i\theta})$ with respect to the measure $dsd\theta$. We will denote by $x_s$ the process
$$
s\in \R_+\to x_s:=X_{\hat g,e^{-s}}(0).
$$

 We have arrived at
 the following useful decomposition of the chaos measure around $z_1=0$:

\begin{lemma}\label{decompmeas}
On the ball $B(0,1)$ we have the following decomposition of the measure $M_\gamma$:
$$
 \int_A\frac{1}{|x|^{\gamma Q}}\,M_\gamma(dx) =c_\gamma \int_0^\infty \int_0^{2\pi}\mathbf{1}_A(e^{-s}e^{i\theta})e^{\gamma x_s}\hat g(e^{-s})\,\mu_{Y}(ds,d\theta)$$
for all $A\subset B(0,1)$ where 
$\mu_{Y}(ds,d\theta)$ is a measure independent of the whole process $(x_s)_{s\geq 0}$. Furthermore, for all $q\in]-\infty;\frac{4}{\gamma^2}[$, we have
\begin{equation}\label{momentbound}
\sup_{a>0}\E\Big[ \Big(\int_a^{a+1}\int_0^{2\pi}e^{\gamma(x_s-x_{a})} \mu_Y(ds,d\theta)\Big)^{q}\Big]<+\infty.
\end{equation}
\end{lemma}

\noindent {\it Proof.} We have for $q< \frac{4}{\gamma^2}$
\begin{align*}
& \E\Big[ \Big(\int_{a}^{a+1}\int_0^{2\pi}e^{\gamma(x_s-x_{a} )} \mu_Y(ds,d\theta)\Big)^{q}\Big]   \\
& \leq (2\pi)^q \E\Big[ e^{q \gamma  \sup_{\sigma \in [0,1]}(x_{a+\sigma}-x_{a} )}  \Big]  
 \E\Big[   \mu_Y( [a,a+1] \times [0,2 \pi] )^q\Big ]  \\
& =  (2\pi)^q \E\Big[e^{q \gamma  \sup_{\sigma \in [0,1]}(x_{a+\sigma}-x_{a} )}  \Big]  
 \E\Big[   \mu_Y( [0,1] \times [0,2 \pi] )^q \Big],  \\
\end{align*} 
by stationarity of  $(s,\theta)\in \R_+^*\times[0,2\pi]\mapsto  Y(e^{-s} e^{i\theta})$. By Lemma \ref{BM} the first exponent is Brownian motion and hence the expectation is bounded uniformly in $a$. From Gaussian multiplicative chaos theory \cite[Theorem 2.12]{review}, we have finiteness of the quantity $\E\Big[   \mu_Y( [0,1] \times [0,2 \pi] )^q \Big] < \infty$,   hence we get \eqref{momentbound}.

%

\qed

It will be useful in the proofs to introduce  for all $a\geq 1$ the stopping times $T_a$ defined by
\begin{equation}\label{stoppingtime}
T_{a}=\inf\{s; \: x_s\geq a-1\},
\end{equation}
and we denote by $\mathcal{G}_{T_a}$ the associated filtration.
We have the following analog of \eqref{momentbound} with stopping times

\begin{lemma}\label{decompmeas11}
For all $q \leq 0$, $n\geq 1$,
  \begin{equation}\label{momentbound1}
\E\Big[ \Big(\int_{T_{n-1}}^{T_{n}}\int_0^{2\pi}e^{\gamma (x_{s} -x_{T_{n-1}})} \mu_Y(ds,d\theta)\Big)^{q}\Big]<\infty.
\end{equation}
\end{lemma}
 
 \proof
 Using the independence of the processes $x_{r}$ and $Y$, Lemma \ref{BM} and stationarity of $Y(s,\theta)$ in $s$ we see that \eqref{momentbound1} is equivalent to proving
  \begin{equation}\label{momentbound11}
\E\Big[ \Big(\int_0^\tau\int_0^{2\pi}e^{\gamma \beta_{s}} \mu_Y(ds,d\theta)\Big)^{q}\Big]<\infty.
\end{equation}
 where $\beta$ is a Brownian motion independent of $Y$ and $\tau=\inf\{s;\beta_s\geq 1\}$.
 We have 
 \begin{align*}
 & \E\Big[ \Big(\int_0^{\tau}\int_0^{2\pi}e^{\gamma \beta_s } \mu_Y(ds,d\theta)\Big)^{q}\Big] \leq
   \E\Big[\mathbf{1}_{\tau\leq 1} \Big(\int_0^{\tau}\int_0^{2\pi}e^{\gamma \beta_s } \mu_Y(ds,d\theta)\Big)^{q}\Big] +
    \E\Big[ \Big(\int_0^{1}\int_0^{2\pi}e^{\gamma \beta_s } \mu_Y(ds,d\theta)\Big)^{2q}\Big]^\hf 
 \end{align*}
 The second term is bounded by Lemma \ref{decompmeas}. The first one equals
 \begin{align*}
 &     \sum_{n \geq 1}  \E\Big[ \mathbf{1}_{1/2^{n+1}  <  \tau \leq 1/{2^n}} \Big(\int_0^{\tau}\int_0^{2\pi}e^{\gamma \beta_s } \mu_Y(ds,d\theta)\Big)^{q}\Big]   \\
& \leq  \sum_{n \geq 1}  \E\Big[ \mathbf{1}_{1/2^{n+1}  <  \tau \leq 1/{2^n}}  \Big(\int_0^{1/2^{n+1}}\int_0^{2\pi}e^{\gamma \beta_r } \mu_Y(dr,d\theta)\Big)^{q}\Big]   \\
& \leq  \sum_{n \geq 1}  \P(1/2^{n+1}  <  \tau \leq 1/{2^n} )^{1/2} \E\Big[  \Big(\int_0^{1/2^{n+1}}\int_0^{2\pi}e^{\gamma \beta_r } \mu_Y(dr,d\theta)\Big)^{2q}\Big]^{1/2}   \\
& \leq  \sum_{n \geq 1}  \P(1/2^{n+1}  <  \tau \leq 1/{2^n} )^{1/2}
 \E\Big[e^{q \gamma  \sup_{\sigma \in [0,2^{-n-1}]}\beta(\sigma)}  \Big]^\hf  
 \E\Big[   \mu_Y( [0,2^{-n-1}] \times [0,2 \pi] )^q\Big]^\hf   \\
& \leq C \sum_{n \geq 1} e^{-c2^n}
 \E\Big[   \mu_Y( [0,2^{-n-1}] \times [0,2 \pi] )^q\Big ]^\hf   \\
& \leq C \sum_{n \geq 1} e^{-c2^n} 
 \E\Big[   \mu_Y( [0,2^{-n-1}] \times [0,2^{-n-1}] )^q \Big]^\hf. 
\end{align*} 
One can find some constant $C>0$ such that the covariance $\E[Y(e^{-s} e^{i \theta}) Y(e^{-s'} e^{i \theta'})   ]$ is bounded by $\ln \frac{1}{|se^{i \theta}-s'e^{i \theta'}|}+C$ hence by Kahane's convexity inequality \cite[Theorem 2.1]{review} one gets the existence of some constant $C>0$ such that
\begin{equation*}
\E\Big[   \mu_Y( [0,2^{-n-1}] \times [0,2^{-n-1}] )^q \Big] \leq C \frac{1}{2^{n \xi(-q)}},
\end{equation*}
 with $\xi(-q)=-(2+\frac{\gamma^2}{2})q - \gamma^2 \frac{q^2}{2}$.
Hence $\sum_{n \geq 1} e^{-c2^n} 
 \E\Big[   \mu_Y( [0,2^{-n-1}] \times [0,2^{-n-1}] )^q ]^\hf < \infty$, which concludes the proof.
\qed

Finally, we will consider a probability measure associated to the martingale $( f_\epsilon^{n} )_{\epsilon \in]0,1]}$ defined by
\begin{equation}\label{def:fn}
 f_\epsilon^{n}  =\mathbf{1}_{ \big\{\min_{u\in[\epsilon,1]}n -X_{\hat{g},u}(0)   \geq 0\big\}}(n -x_{\ln \frac{1}{\epsilon}}  ).
\end{equation}
The martingale property of  $( f_\epsilon^{n} )_{\epsilon \in]0,1]}$ is classical and results from Lemma \ref{BM} as well as the stopping time theorem.
We can define for each $\epsilon\in]0,1]$ a probability measure on $\mathcal{F}_\epsilon$ by
$$\Theta^{n}_\epsilon =\frac{1}{\E[ f_\epsilon^{n}  ]} f_\epsilon^{n}  \,d\P,$$
where one has the following bound $\E[ f_\epsilon^{n }  ]=\E[ f_1^{n }  ] \leq  n+C$ for some constant $C$.
Because of Lemma \ref{BM} and the martingale property of the family $( f_\delta^{n } )_{\delta \in]0,1]}$, it is plain to check that these probability measures are compatible in the sense that, for $\epsilon'<\epsilon$
\begin{equation}\label{def:thetaeps}
\Theta^{n }_{\epsilon'} |_{\mathcal{F}_\epsilon}=\Theta_\epsilon^{n }.
\end{equation}
By Caratheodory's extension theorem we can find a probability measure $\Theta^{n }$ on $\mathcal{F}_\infty$ such that for all $\epsilon \in]0,1]$
\begin{equation}\label{def:theta}
\Theta^{n }|_{\mathcal{F}_\epsilon}=\Theta^{n }_\epsilon.
\end{equation}
We denote by $\E^{\Theta^{n }}$ the corresponding expectation.

Recall the following explicit law of the Brownian motion conditioned to stay positive 

\begin{lemma}\label{Bessel}
Under the probability measure $\Theta^{n }$, the process  
$$ t\mapsto   n  -x_t  $$
  evolves  as a $3d$-Bessel process starting from $ n-x_0 $ where $x_0$ is distributed like $X_{\hat{g},1}$ (under $\P$) conditioned to be less or equal to $n $.
\end{lemma} 

%

We will sometimes use the following classical representation: under  $\Theta^{n}$, the  process $ t\mapsto   n  -x_t $ is distributed like $| n-x_0+B_t|$ where $B_t$ is a standard 3d Brownian motion starting from $0$ (here, we identify $ n-x_0$ with $( n-x_0)(1,0,0)$).

\section{Construction of the derivative Q-vertex}
In this section, we prove the claims in  Lemma \ref{crucial1} concerning $\tilde A_\epsilon$.
\medskip

\noindent {\it Proof of \eqref{est4}}. From \eqref{Bep} we get
\begin{align}
|B_\epsilon(F,n)|\leq &C \int_\R e^{\sigma c } \E \Big[\mathbf{1}_{M_{n,\epsilon}}  |n+c|  e^{ - \mu e^{\gamma c} \int_{D_\epsilon}e^{\gamma H_\epsilon} \,dM_\gamma } \Big]\,dc.
\end{align}
Recall  that $|z_i|>1$ for $i\geq 2$. Then, recalling \eqref{def:Hepsilon} and \eqref{hatGformula} we get for $\epsilon\leq |z|\leq 1$
\begin{equation}\label{minH}
e^{\gamma H_\epsilon(z)}\geq C|z|^{-\gamma Q}.
\end{equation}
By Lemma \ref{decompmeas}, we get

\begin{align}\label{Bepsest}
| B_\epsilon(F,n) |\leq 
&C \int_\R e^{\sigma c }  \,\E\Big[\mathbf{1}_{M_{n,\epsilon}} |n+c| 
e^{-C \mu e^{\gamma c} \int_{0}^{\ln\frac{1}{\epsilon}} e^{\gamma  x_{r}}\, \mu_{Y}(dr)  } \Big]\,dc ,
\end{align}
where $\mu_{Y}(dr)$ is the measure defined by $\mu_{Y}(dr)=\int_0^{2\pi}\mu_{Y}(dr,d\theta)$.  

Below, we want to show that the integral in the exponential term above carries a big amount of mass, and we will look for this mass at some place where the process $r\mapsto x_{r}$  takes on values close to its maximum, which is between $n-1$ and $n$ on the set $M_{n,\epsilon}$. To locate this place, we use the stopping times $T_{n-1}$ and $T_{n}$ defined by \eqref{stoppingtime} which are finite and belong to $[0,\ln\frac{1}{\epsilon}]$ on $M_{n,\epsilon}$. We deduce

\begin{align}
 | B_\epsilon(F,n)|
\leq &C \int_\R e^{\sigma c } \,\E\Big[\mathbf{1}_{M_{n,\epsilon}}   |n+c|     e^{- \mu e^{\gamma c} Ce^{\gamma(n-1)} I_n } \Big]\,dc \nonumber
\end{align}
where we have set 
\begin{equation}\label{Indef}
I_n=\int_{T_{n-1}}^{T_n } e^{\gamma  (x_{r}-x_{T_{n-1}})}\,\mu_{Y}(dr).
\end{equation}
 By making the change of variables $y=e^{\gamma(c+n)}I_n$, we get 
 
\begin{align}
   B_\epsilon(F,n)\leq &C  e^{-n\sigma } \int_0^{\infty} y^{\frac{\sigma }{\gamma}-1  } (1+|\ln y|) e^{-\mu C e^{-\gamma}y}\,dy\, \E\Big[\mathbf{1}_{M_{n,\epsilon}}    (1+|\ln I_n|)I_n^{-{\frac{\sigma }{\gamma}}}\Big] . \nonumber
\end{align}
Then we bound  

\begin{equation*}
 \E\Big[\mathbf{1}_{M_{n,\epsilon}}    (1+|\ln I_n|)I_n^{-{\frac{\sigma }{\gamma}}}\Big]
 \leq \P (M_{n,\epsilon})^{1/2}  \E\Big[\mathbf{1}_{M_{n,\epsilon}}    (1+|\ln I_n|)^2I_n^{-{\frac{2\sigma }{\gamma}}}\Big]^\hf
.
\end{equation*}
Hence, by Lemma \ref{decompmeas11} we conclude 
\begin{align}
 B(F;\epsilon,n)\leq &C  e^{-n\sigma } \P (M_{n,\epsilon})^{1/2} . \nonumber
\end{align}

 The claim \eqref{est4} then follows  by the dominated convergence theorem since for each fixed $n$, the probability $\P (M_{n,\epsilon})$ goes to $0$ as $\epsilon$ goes to $0$ (see Lemma \ref{standard}).
\qed

\subsubsection{Proof of  \eqref{est3}.}
Proceeding as in the  proof of  \eqref{est4} we get 
\begin{equation*}  
 \tilde A_\epsilon(1,n)\leq C e^{-n\sigma }\E\Big[\mathbf{1}_{M_{n,\epsilon}} ( n-x_{\ln \frac{1}{\epsilon}})    I_n^{-{\frac{\sigma }{\gamma}}} \Big],  
\end{equation*}
where $I_n$ is as in \eqref{Indef}. Now, we have
\begin{align*}
 \E[  \mathbf{1}_{M_{n,\epsilon}} ( n-x_{\ln \frac{1}{\epsilon}}) | \mathcal{F}_{T_n} \vee \sigma(Y) ]   
 & = \E[  \mathbf{1}_{\min_{s \in [0,\ln \frac{1}{\epsilon}]}( n-x_s) \geq 0} ( n-x_{\ln \frac{1}{\epsilon}}) | \mathcal{F}_{T_n} \vee \sigma(Y)  ]  \mathbf{1}_{T_n \leq \ln \frac{1}{\epsilon}}  \\
&=  \mathbf{1}_{\min_{s \in [0,T_n]}( n-x_s) \geq 0}
 ( n-x_{T_n}) \mathbf{1}_{T_n \leq \ln \frac{1}{\epsilon}}  \\
& =   \mathbf{1}_{\min_{s \in [0,T_n]}( n-x_s) \geq 0}
 \mathbf{1}_{T_n \leq \ln \frac{1}{\epsilon}} 
\end{align*}
so that
\begin{equation*}
 \tilde A_\epsilon(1,n) \leq Ce^{-n\sigma } \E\, I_n^{-{\frac{\sigma }{\gamma}}}
 \leq  C e^{-n\sigma }.   
\end{equation*}
from which the estimate \eqref{est3} follows.

\subsubsection{Proof of  the first part of \eqref{two limits}, i.e. the existence of $\lim_{\epsilon\to 0} \tilde A_{\epsilon}(F,n)$.}

Now, we need to establish the existence and non triviality of the limit of $\tilde A_{\epsilon}(F,n)$, i.e. one part of \eqref{two limits}.
Since $H_\epsilon$ converges in $H^{-1}(\S^2)$ towards $H$, it suffices to study the convergence and non triviality of the limit for $F=1$ and fixed $n$. We claim that this will result from the convergence in probability of the quantity $\int_{D_\epsilon}e^{\gamma H_\epsilon }\,dM_\gamma$ under the probability measure $\Theta^{n}$ towards a non trivial limit. To see this,  make the change of variables $y=e^{\gamma c}\int_{D_\epsilon}e^{\gamma H_\epsilon }\,dM_\gamma$ to get

\begin{align*}\label{dom}
 \int e^{\sigma c}\E\Big[  &  \mathbf{1}_{M_{n,\epsilon}}  (n-x_{\ln \frac{1}{\epsilon}})  e^{ - \mu e^{\gamma c}\int_{D_\epsilon}e^{\gamma H_\epsilon }\,dM_\gamma } \Big]\,dc\\
 =&{\gamma}^{-1}\E[ f_1^{n}  ] \int_{0}^{\infty}y^{\frac{\sigma }{\gamma}-1}e^{-\mu y}\,dy\times  \E^{\Theta^{n}}\Big[ \mathbf{1}_{M_{n,\epsilon}}  \Big(\int_{D_\epsilon}e^{\gamma H_\epsilon }\,dM_\gamma\Big)^{-\frac{\sigma }{\gamma}} \Big].
%
\end{align*}
Under the probability measure $\Theta^{n}$ the process $t \mapsto (n-x_t)$ is a 3d Bessel process hence $\min_{s \in [0,\ln \frac{1}{\epsilon}]}( n-x_s)$  converges almost surely to a finite random variable as $\epsilon$ goes to $0$ and therefore $\mathbf{1}_{M_{n,\epsilon}}$ converges to $\mathbf{1}_{\max_{s \in [0,\infty]}(x_s) \in [n-1,n]}$.

Take any non empty closed ball $B$ of $\R^2$ containing no insertions $z_i$. Then 
$\sup_\epsilon H_\epsilon$ is bounded in $B$ and thus
$$ 
\Big(\int_{D_\epsilon}e^{\gamma H_\epsilon }\,dM_\gamma\Big)^{-\frac{\sigma }{\gamma}} \leq  CM_\gamma(B)^{-\frac{\sigma }{\gamma}} .
$$
Let $\delta>0$ be such that $B\subset B(0,\delta)^c$. 
Then 
$$
\E^{\Theta^{n}}\Big[ M_\gamma(B)^{-\frac{\sigma }{\gamma}} \Big] \leq C (n+1)^{-1}\E[f^{n}_\delta M_\gamma(B) ^{-\frac{\sigma }{\gamma}} ]\leq C (n+1)^{-1}\E[(f^{n}_\delta)^2]^{\hf}\E\Big[ M_\gamma(B)^{-\frac{2\sigma }{\gamma}} ]^{\hf}.
 $$
Because GMC admits moments of negative order \cite[theorem 2.12]{review}, the last expectation is finite. Hence the dominated convergence theorem entails that to prove our claim it is enough to establish the convergence in probability of the quantity 
$\int_{D_\epsilon}e^{\gamma H_\epsilon }\,dM_\gamma$ under the probability measure $\Theta^{n}$ towards a non trivial limit. Because $M_\gamma$ is a positive measure and because of the bound   \eqref{minH},  this is clearly equivalent to the finiteness under $\Theta^{n}$ of the quantity $\int_{\R^2}e^{\gamma H }\,dM_\gamma$.  Outside of the ball $B(0,1)$, the finiteness results from the fact that $  \int_{D_1}e^{\gamma H}\,dM_\gamma<\infty$ under $\P$ (see see  \cite[proof of Th. 3.2]{DKRV}), and the absolute continuity of  $\Theta^{n}$ with respect to $\P$ when restricted to $\mathcal{F}_1$. The main point  is thus to analyze the integrability inside the ball $B(0,1)$. It is  clearly enough to show
\begin{equation}\label{sing}
 \int_{B(0,1)}\frac{1}{|x|^{\gamma Q}}\,dM_\gamma <\infty, \quad \text{a.s. under }\Theta^{n }.
\end{equation}
This follows  from  the following Lemma \ref{integM}:
     

\qed

 \begin{lemma}\label{integM}
The measure $M_\gamma$ satisfies 
\begin{equation}
\E^{\Theta^{n} }\Big[ \int_{B(0,1)}\frac{1}{|x|^{\gamma Q}}M_\gamma(dx)\Big]<\infty.
\end{equation}
\end{lemma}

\noindent{\it Proof.} 
Under the measure $\Theta^{n }$, the process $ t\mapsto   n  -x_t $ is distributed like $ |n-x_0+B_t|$ where $B_t$ is a standard $3$ dimensional Brownian motion (here, we identify $ n-x_0$ with $(n-x_0) (1,0,0)$). We suppose the Brownian motion lives on the same probability space. Then, if $N$ denotes a standard 3d Gaussian variable (under some expectation we will also denote $\E$), we have
\begin{align*}
 \E^{\Theta^{n } }\Big[ \int_{B(0,1)}\frac{1}{|x|^{\gamma Q}}M_\gamma(dx)\Big]  & \leq C \E^{\Theta^{n} }\Big[ \int_0^\infty \int_0^{2\pi}e^{\gamma x_{r}}\,\mu_{Y}(dr,d\theta)  \Big ]   \\
& = C  e^{\gamma n} \E^{\Theta^{n } }\Big[ \int_0^\infty \int_0^{2\pi}e^{-\gamma |n-x_0+B_r| }\,\mu_{Y}(dr,d\theta)  \Big ]     \\
& = C  e^{\gamma n} \E^{\Theta^{n } }\Big[ \int_0^\infty e^{-\gamma |n-x_0+B_r| }\, dr \Big ]     \\
& \leq C e^{2 \gamma  n} \E\,[ e^{\gamma |x_0|} ] \E^{\Theta^{n } }\Big[ \int_0^\infty e^{-\gamma |B_r| }\, dr \Big ]     \\
& = C e^{2 \gamma n} \E\Big[ \int_0^\infty e^{-\gamma \sqrt{r} |N| }\, dr \Big ],     \\
&=   C  e^{2 \gamma  n} ( \int_0^\infty  e^{-\gamma \sqrt{r} }\, dr
 )  \E \Big [\frac{1}{|N|^2}  \Big] < \infty.\qed
\end{align*}

\section{Renormalization of the $Q$-puncture vertex operators}\label{sec:renorm}


\subsection{Proof of  \eqref{est2}} 
Using \eqref{minH} and proceeding as for \eqref{Bepsest} we  get
\begin{align}
  {A}_\epsilon(1,n) 
\leq &C\int_\R e^{\sigma c }   \,\E\Big[\mathbf{1}_{M_{n,\epsilon}}    \exp\Big( - \mu e^{\gamma c} C \int_{0}^{\ln\frac{1}{\epsilon}} e^{\gamma  x_r}\,\mu_{Y}(dr)  \Big) \Big]\,dc \nonumber
\end{align}
 The stopping time 
$T_{n}=\inf\{s;\: x_s\geq n-1\}$
 is finite and belongs to $[0,\ln\frac{1}{\epsilon}]$ on $M_{n,\epsilon}$. We deduce that
 
\begin{align}
 {A}_\epsilon(1,n) 
\leq &C\int_\R e^{\sigma c }     \,\E\Big[\mathbf{1}_{M_{n,\epsilon}\cap \{T_n<\ln\frac{1}{\epsilon}-1\}}    \exp\Big( - \mu e^{\gamma c} Ce^{\gamma(n-1)} I(T_n)  \Big) \Big]\,dc \\
&+C\int_\R e^{\sigma c }     \,\E\Big[\mathbf{1}_{M_{n,\epsilon}\cap \{T_n\geq \ln\frac{1}{\epsilon}-1\}}    \exp\Big( - \mu e^{\gamma c} Ce^{\gamma x_{\ln\frac{1}{\epsilon}-1} I(\ln\frac{1}{\epsilon}-1)}  \Big) \Big]\,dc \nonumber\\
&=: a_\epsilon(n)+ b_\epsilon(n)
\end{align}
where we have set $$I(z)=\int_{z}^{z+1} e^{\gamma  (x_r-x_z)}\,\mu_{Y}(dr).$$ 
We will show that there exists a constant $C>0$ such that for all $n$ 
\begin{equation}
 (\ln\frac{_1}{^\epsilon})^{\hf}a_\epsilon(n), \  (\ln\frac{_1}{^\epsilon})^{\hf}b_\epsilon(n)
\leq C n e^{-\sigma n},
\end{equation}
which is enough to complete the proof of  \eqref{est2}. 

\medskip
We begin with $ a_\epsilon(n) $.
By making the change of variables $y=e^{\gamma(c+n)}I(T_n)$, we get
\begin{align}
  a_\epsilon(n)\leq &C\,  e^{-n\sigma } \int_0^{\infty} y^{\frac{\sigma }{\gamma}-1  }  e^{-\mu C e^{-\gamma}y}\,dy \E\Big[\mathbf{1}_{M_{n,\epsilon}\cap \{T_n<\ln\frac{1}{\epsilon}-1\}}    I(T_n)^{-{\frac{\sigma }{\gamma}}} \Big] . \nonumber
\end{align}
It suffices to estimate the last expectation. Obviously, we have

\begin{align}\label{tam}
 \E\Big[&\mathbf{1}_{M_{n,\epsilon}\cap \{T_n<\ln\frac{1}{\epsilon}-1\}}  I(T_n)^{-{\frac{\sigma }{\gamma}}}  \Big] \\
\leq & \E\Big[\mathbf{1}_{\{\min_{u\in [0,T_n]}n-x_u\geq 0\}} \mathbf{1}_{\{\min_{u\in [T_n+1,\ln\frac{1}{\epsilon}]}(n-x_{T_n+1})-(x_u-x_{T_n+1})\geq 0\}}  \frac{\mathbf{1}_{\{T_n+1< \ln\frac{1}{\epsilon}\}} }{I(T_n)^{\frac{\sigma }{\gamma}}} \Big].\nonumber
 \end{align}
By conditioning on the  the sigma algebra  $\mathcal{H}_{T_n}$ generated by 
$\{x_{r},r\leq T_n\}$, $\{x_{r}-x_{T_n+1},  r\geq  T_n+1\}$ and $\{x_{T_n+1}-n\}$, we see that we have to estimate the quantity
$$\E[I(a)^{-\frac{\sigma }{\gamma}}|x_{a+1}-x_{a}].$$
We claim

\begin{lemma}\label{bridge}
There exists a constant $C$ (independent of any relevant quantity) such that for all $a>0$
$$\E[I(a)^{-\frac{\sigma }{\gamma}}|x_{a+1}-x_a]\leq C\big(e^{-\sigma (x_{a+1}-x_a)}+1\big).$$
\end{lemma}

The proof of this lemma is given just below. Admitting it for a while and given the fact that the random variable $x_{T_n+1}-n$ is a standard Gaussian random variable, the conditioning on $\mathcal{H}_{T_n }$ of the expectation \eqref{tam} thus gives

\begin{align*}
 \E\Big[&\mathbf{1}_{M_{n,\epsilon}\cap \{T_n<\ln\frac{1}{\epsilon}-1\}}    \frac{1}{I(T_n)^{\frac{\sigma }{\gamma}}} \Big] \\
\leq & C\int_{\R}\E\Big[\mathbf{1}_{\{\min_{u\in [0,T_n]}n-x_{u}\geq 0\}} \mathbf{1}_{\{\min_{u\in [T_n+1,\ln\frac{1}{\epsilon}]}-y-(x_{u}-x_{T_n+1})\geq 0\}} \Big]  \big(e^{-\sigma (y+1)}+1\big) e^{-y^2/2}\,dy
 \end{align*}
 To estimate the expectation in the integral, use the strong Markov property of the Brownian motion to write
 
 \begin{align*}
\E\Big[&\mathbf{1}_{\{\min_{u\in [0,T_n]}n-x_{u}\geq 0\}} \mathbf{1}_{\{\min_{u\in [T_n+1,\ln\frac{1}{\epsilon}]}-y-(x_{u}-x_{T_n+1})\geq 0\}} \Big]\\
=&\E\Big[\mathbf{1}_{\{\min_{u\in [0,T_n]}n-x_{u}\geq 0\}} \mathbf{1}_{
\{\min_{u\in [T_n,\ln\frac{1}{\epsilon}-1]}-y-(x_{u}-x_{T_n})\geq 0\}} \Big]\\
\leq & \E\Big[\mathbf{1}_{\{\min_{u\in [0,\ln\frac{1}{\epsilon}-1]}n+\max (0,-y)-x_{u}\geq 0\}}  \Big]
\leq  (\frac{_2}{^\pi})^\hf \frac{n+\max (0,-y)}{(\ln\frac{1}{\epsilon}-1)^\hf},
 \end{align*}
where in the last inequality we have used Lemma \ref{standard}.
 We deduce
 $$
 \E\Big[\mathbf{1}_{M_{n,\epsilon}\cap \{T_n<\ln\frac{_1}{^\epsilon}-1\}}    \frac{1}{I(T_n)^{\frac{\sigma }{\gamma}}} \Big] 
\leq C(\ln\frac{_1}{^\epsilon})^{-\hf}n  e^{-n\sigma }
 $$

All in all, we have obtained
 $$ \sup_{\epsilon\in]0,1]}(\ln\frac{_1}{^\epsilon})^{\hf} a_\epsilon(n)\leq Cn  e^{-n\sigma } ,$$
 which proves the claim. The same argument holds for $b_\epsilon(n)$.\qed

 \medskip

\noindent {\it Proof of Lemma \ref{bridge}.} Notice that the joint law of $\Big((x_r-x_a)_{r\in[a,a+1]},x_{a+1}-x_a\Big)$ is that of $\big((B_u-B_a)_{u\in[a,a+1] },B_{a+1}-B_a)$ where $B$ is a standard Brownian motion starting from $0$ (independent of $Y$). Hence the law of $I(a)$ conditionally on $x_{a+1}-x_a=x$ is given by
$$\int_{0}^{1}e^{\gamma {\rm Bridge}^{0,x}_r}\,\mu
_Y(dr)$$ where $({\rm Bridge}^{0,x}_r)_{r\leq 1}$ is a Brownian bridge between $0$ et $x$ with lifetime $1$. Hence it has the law of $r\mapsto B_r-rB_1+ux$. By convexity of the mapping $x\mapsto x^{-q}$ for $q>0$ and the fact that the covariance kernel of the Brownian Bridge and the Brownian motion are comparable up to fixed constant, we can apply Kahane's  inequality \cite{cf:Kah} to get that \begin{align*}
\E[I(a)^{-\frac{\sigma }{\gamma}}|x_{a+1}-x_a=x]\leq C\E\Big[\Big(\int_{a}^{a+1} e^{\gamma ( B_r-B_a)+(r-a)x}\,\mu_{Y}(dr)\Big)^{-\frac{\sigma }{\gamma}}\Big]
\end{align*}
From  Lemma \ref{decompmeas} and the fact that $e^{(r-a)x}\geq e^x\wedge 1$ for $r\in[a,a+1]$, this quantity is less than
$$\E[I(a)^{-\frac{\sigma }{\gamma}}|x_{a+1}-x_a=x]\leq C\big(e^{-\sigma x}\vee 1\big).$$
This proves the claim.\qed

 \subsection{Proof of  \eqref{two limits}.} 
 First notice that
\begin{align*} 
 \sum_{n=0}^N{A}_\epsilon(1,n)=&   \int_\R e^{\sigma c }  \,\E\Big[\mathbf{1}_{B_{N,\epsilon}} 
    \exp\Big( - \mu e^{\gamma c} \int_{D_\epsilon}e^{\gamma H_\epsilon }\,dM_\gamma \Big) \Big]\,dc
\end{align*}
where
$$B_{N,\epsilon}=\{\min_{u\in[\epsilon,1]}  N-x_u\geq 0\}.$$
Let us denote by $Z_\epsilon$ the measure $e^{\gamma H_\epsilon }\,dM_\gamma$ and define 
$$
s_\epsilon:=(\ln\frac{_1}{^\epsilon})^{\frac{_1}{^6}},\ \ h_\epsilon:=e^{-s_\epsilon}.$$ 
Now we prove the  upper bound. We have
\begin{align*} 
 \sum_{n=0}^N{A}_\epsilon(1,n)\leq& \int_\R e^{\sigma c }  \,\E\Big[\mathbf{1}_{B_{N,\epsilon}} 
    \exp\Big( - \mu e^{\gamma c} Z_\epsilon(D_{h_\epsilon})  \Big) \Big]\,dc\\
    =&  \int_\R e^{\sigma c }  \,\E\Big[\E\big[\mathbf{1}_{B_{N,\epsilon}} |\mathcal{F}_{s_\epsilon}
  \big]  \exp\Big( - \mu e^{\gamma c}Z_\epsilon(D_{h_\epsilon})\Big) \Big]\,dc.
\end{align*}
Using  the standard estimate $\E\big[\mathbf{1}_{B_{N,\epsilon}} |\mathcal{F}_{s_\epsilon}  \big] \leq \sqrt{2/\pi}\frac{ N-x_{s_\epsilon}}{\sqrt{\ln\frac{_1}{^\epsilon}-s_\epsilon}}$ (see Lemma \ref{standard}) we deduce
$$ \limsup_{\epsilon\to 0} \sum_{n=0}^N (\ln\frac{_1}{^\epsilon})^\hf{A}_\epsilon(1,n)\leq   \sqrt{2/\pi}  \sum_{n=0}^N\tilde A(1,n).$$
which completes the upper bound.

\medskip  Let us now investigate the lower bound. We denote by $C(\epsilon)$ the annulus $\{x:\epsilon\leq |x|\leq h_\epsilon\}$ and by $I_\epsilon$ the set 
$$
I_\epsilon=\{\min_{u\in[0,s_\epsilon]}( N-x_u)\geq s_\epsilon^\theta\} 
$$
where $\theta\in]0,1/2[$. We have

\begin{align*}
\sum_{n=0}^N{A}_\epsilon(1,n)\geq   \int_\R e^{\sigma c }  \,\E\Big[\mathbf{1}_{B_{N,\epsilon}}  \mathbf{1}_{I_{\epsilon}} 
    \exp\Big( - \mu e^{\gamma c}Z_\epsilon(D_{h_\epsilon})- \mu e^{\gamma c}Z_\epsilon (C(\epsilon)\Big) \Big]\,dc. 
\end{align*}
Using  $e^{-u}\geq 1-u^\hf$, we deduce
\begin{align}
\sum_{n=0}^N{A}_\epsilon(1,n)\geq &  \int_\R e^{\sigma c }  \,\E\Big[\mathbf{1}_{B_{N,\epsilon}}\mathbf{1}_{I_{\epsilon}} 
    e^{ - \mu e^{\gamma c}Z_\epsilon(D_{h_\epsilon}) } \Big(1-\mu^\hf e^{\hf\gamma c}Z_\epsilon (C(\epsilon))^\hf\Big)\Big]\,dc\nonumber \\
  =&  \int_\R e^{\sigma c }  \,\E\Big[\mathbf{1}_{B_{N,\epsilon}}  
    e^{ - \mu e^{\gamma c}Z_\epsilon(D_{h_\epsilon})  } \Big]\,dc
     -\int_\R e^{\sigma c }  \,\E\Big[\mathbf{1}_{B_{N,\epsilon}}  \mathbf{1}_{I_{\epsilon}^c} 
    e^{- \mu e^{\gamma c}Z_\epsilon(D_{h_\epsilon}) } \Big]\,dc\nonumber \\
    &-  \mu^\hf\int_\R e^{c\big( \hf\gamma+\sigma \big) }  \,\E\Big[\mathbf{1}_{B_{N,\epsilon}} \mathbf{1}_{I_{\epsilon}} e^{- \mu e^{\gamma c} Z_\epsilon(D_{h_\epsilon}) } Z_\epsilon(C(\epsilon))^\hf \Big]\,dc\nonumber\\    
 =:
 &{B}_1(N,\epsilon)-{B}_2(N,\epsilon)-{B}_3(N,\epsilon).\label{sumB}
\end{align}
We now estimate the above three terms. 

\medskip

We start with ${B}_1(N,\epsilon)$.
We have 
\begin{align*}
\E\big[\mathbf{1}_{B_{N,\epsilon}}  |\mathcal{F}_{h_\epsilon}  \big]=&\mathbf{1}_{B_{N,h_\epsilon}}(\frac{_2}{^\pi})^\hf\int_0^{\frac{ N-x_{s_\epsilon}}{\sqrt{\ln\frac{_1}{^\epsilon}-s_\epsilon}}}e^{-\frac{u^2}{2}}\,du
\\
\geq & \mathbf{1}_{B_{N,h_\epsilon}}\mathbf{1}_{\{ N-x_{s_\epsilon}\leq 
(\ln\frac{_1}{^\epsilon}-s_\epsilon)^{\frac{1}{4}}\}}(\frac{_2}{^\pi})^\hf
\int_0^{\frac{ N-x_{s_\epsilon}}{\sqrt{\ln\frac{_1}{^\epsilon}-s_\epsilon}}}e^{-\frac{u^2}{2}}\,du
\\
\geq &  \mathbf{1}_{B_{N,h_\epsilon}}\frac{ N-x_{s_\epsilon}}{\sqrt{\ln\frac{_1}{^\epsilon}-s_\epsilon}}\mathbf{1}_{\{ N-x_{s_\epsilon}\leq (\ln\frac{_1}{^\epsilon}-s_\epsilon)^{\frac{1}{4}}\}} (\frac{_2}{^\pi})^\hf e^{-\hf (\ln\frac{_1}{^\epsilon}-s_\epsilon)^{-\hf}} .
\end{align*}
 Plugging this relation into ${B}_1(N,\epsilon)$ we deduce
 
 \begin{align*}
 {B}_1(N,\epsilon) &\geq (\frac{_2}{^\pi})^\hf e^{-\hf (\ln\frac{_1}{^\epsilon}-s_\epsilon)^{-\hf}}
 (\ln\frac{_1}{^\epsilon}-s_\epsilon)^{-\hf}\big(
 \int_\R e^{\sigma c }  \,\E\Big[ \mathbf{1}_{B_{N,\epsilon}}( N-x_{s_\epsilon})  \exp\Big( - \mu e^{\gamma c} Z_\epsilon(D_{h_\epsilon}) \Big) \Big]\,dc
 \\
 &-\int_\R e^{\sigma c }  \,\E\Big[\mathbf{1}_{\{ N-x_{s_\epsilon}>(\ln\frac{_1}{^\epsilon}-s_\epsilon)^{\frac{1}{4}}\}}  \mathbf{1}_{B_{N,h_\epsilon}}( N-x_{s_\epsilon})  e^{ - \mu e^{\gamma c} Z_\epsilon(D_{h_\epsilon})}\Big]\,dc
 \\
 &=:\Delta_1(\epsilon)+\Delta_2(\epsilon).
\end{align*}
 It is clear that 
 $$
 \lim_{\epsilon\to 0}(\ln\frac{_1}{^\epsilon})^\hf\Delta_1(\epsilon)=\lim_{\epsilon\to 0}\sum_{n=0}^N\tilde A_\epsilon(1,n)=\sum_{n=0}^N\tilde A(1,n).
 $$
  It remains to treat  $\Delta_2(\epsilon)$. By making the change of variables $y=e^{\gamma c}Z_\epsilon(D_{h_\epsilon})$, we get
 $$(\ln\frac{_1}{^\epsilon})^\hf\Delta_2(\epsilon)\leq C\E^{\Theta^{N }}\Big[\mathbf{1}_{\{ N-x_{s_\epsilon}> (\ln\frac{_1}{^\epsilon}-s_\epsilon)^{\frac{_1}{^4}}\}} Z_\epsilon(D_{h_\epsilon})^{-\frac{\sigma }{\gamma}}\Big].$$
 Now we will use the fact that under $\Theta^{N } $ the event in the above expectation is very unlikely. Using the elementary inequality $ab\leq a^2/2+b^2/2$ we get
  \begin{align*}
  \E^{\Theta^{N }}&\Big[\mathbf{1}_{\{ N-x_{s_\epsilon}> (\ln\frac{_1}{^\epsilon}-s_\epsilon)^{\frac{_1}{^4}}\}}  Z_\epsilon(D_{h_\epsilon})^{-\frac{\sigma }{\gamma}}\Big]\\
  \leq & (\ln\frac{_1}{^\epsilon})^{ \kappa} \E^{\Theta^{N }}\Big[\mathbf{1}_{\{ N-x_{s_\epsilon}> (\ln\frac{_1}{^\epsilon}-s_\epsilon)^{\frac{_1}{^4}}\}}  \Big]+ (\ln\frac{_1}{^\epsilon})^{-\kappa}\E^{\Theta^{N }}\Big[ Z_\epsilon(D_{1})^{-2\frac{\sigma }{\gamma}}\Big].
 \end{align*}
 Using the fact that a Gaussian Multiplicative Chaos has negative moments of all orders on all open balls, the expectation in the second term  in the above expression is easily seen to be bounded uniformly in $\epsilon$. Hence, the second term tends  to $0$ as $\epsilon\to 0$. 
 Concerning the first term, recall Lemma \ref{Bessel} and the estimate, for a $3d$-Bessel process $\beta_t$ and $u>x$
 $$\P_x(\beta_t>u)=\P_{x/t^\hf}(\beta_1>u/t^\hf)\leq C\frac{{t}^\hf}{u-x}\wedge 1.$$
 Therefore 
  \begin{align*}
   \Theta^{N }\Big( N-x_{s_\epsilon}> (\ln\frac{_1}{^\epsilon}-s_\epsilon)^{\frac{_1}{^4}}\Big)
   \leq & C\E[\frac{s_\epsilon^{1/2}}{|\ln\frac{_1}{^\epsilon}-s_\epsilon
 -x_0|}\wedge 1]
   \\
   \leq & 2C\E[\frac{s_\epsilon^{1/2}}{|\ln\frac{_1}{^\epsilon}-s_\epsilon
 -x_0|}]+\P(x_0>\hf(\ln\frac{_1}{^\epsilon}-s_\epsilon)^{\frac{_1}{^4}})\\
   \leq &2C\ln (\frac{_1}{^\epsilon})^{-\frac{_1}{^6}}+C(\ln\frac{_1}{^\epsilon})^{\frac{_1}{^4}}e^{-\hf (\ln\frac{_1}{^\epsilon})^{\frac{_1}{^2}}}.
  \end{align*}
 Hence, choosing $\kappa<1/6$ leads to $\lim_{\epsilon\to 0}(\ln\frac{_1}{^\epsilon})^\hf\Delta_2(\epsilon)=0$. Thus
\begin{equation}\label{B1}
\liminf_{\epsilon\to 0} (\ln\frac{_1}{^\epsilon})^\hf{B}_1(N,\epsilon)\geq \sum_{n=0}^N\tilde A(1,n).
\end{equation}

\medskip Now we treat $ {B}_3(N,\epsilon)$. To this purpose, we use  first the change of variables $y=e^{\gamma c}Z_\epsilon(D_{h_\epsilon})$ to get
\begin{align*}
{B}_3(N,\epsilon)\leq & C\E [\mathbf{1}_{B_{N,\epsilon}} \mathbf{1}_{I_\epsilon}Z_\epsilon(D_{h_\epsilon})^{-\frac{\hf\gamma+\sigma }{\gamma}}
Z_\epsilon(C(\epsilon))^\hf]\\
=&C \E \Big[\mathbf{1}_{B_{N,\epsilon}} \mathbf{1}_{I_\epsilon}\E[Z_\epsilon(D_{h_\epsilon})^{-\frac{\hf\gamma+\sigma }{\gamma}}
Z_\epsilon(C(\epsilon))^\hf|(x_s)_{s<\infty}]\Big]\\
\leq&C \E \Big[\mathbf{1}_{B_{N,\epsilon}} \mathbf{1}_{I_\epsilon}
\E[Z_\epsilon(D_{h_\epsilon})^{-\frac{\gamma+2\sigma }{\gamma}}|(x_s)_{s<\infty}
]^\hf
\E[
Z_\epsilon(C(\epsilon))|(x_s)_{s<\infty}]^\hf\Big]\\
=&C \E \Big[\mathbf{1}_{B_{N,\epsilon}} \mathbf{1}_{I_\epsilon}\E[Z_\epsilon(D_{h_\epsilon})^{-\frac{\gamma+2\sigma }{\gamma}}|(x_s)_{s<\infty}
]^\hf
 \big(\int_{s_\epsilon}^{\ln\frac{1}{\epsilon}}e^{\gamma x_u}\,du\big)^\hf\Big]
\\
\leq&C \E[Z_\epsilon(D_{h_\epsilon})^{-\frac{\gamma+2\sigma }{\gamma}}]^\hf
\E \Big[\mathbf{1}_{B_{N,\epsilon}} \mathbf{1}_{I_\epsilon}
\int_{s_\epsilon}^{\ln\frac{1}{\epsilon}}e^{\gamma x_u}\,du\Big]^\hf
\\
=&C \E \Big[\mathbf{1}_{B_{N,\epsilon}} \mathbf{1}_{I_\epsilon}\int_{s_\epsilon}^{\ln\frac{1}{\epsilon}}e^{\gamma x_u}\,du\Big]^\hf
  \end{align*}
On the set $I_\epsilon$, we have the estimate
$$
\int_{s_\epsilon}^{\ln\frac{_1}{^\epsilon}}e^{\gamma x_u}\,du \leq C\ln\frac{_1}{^\epsilon} e^{-\gamma s_\epsilon^\theta} =C\ln\frac{_1}{^\epsilon} e^{-\gamma(\ln\frac{_1}{^\epsilon})^{\theta/6}}
$$
which implies
\begin{equation}\label{B3limit}
\lim_{\epsilon\to 0} (\ln\frac{_1}{^\epsilon})^\hf{B}_3(N,\epsilon)=0.
\end{equation}

\medskip 
 Finally we focus on $ {B}_2(N,\epsilon)$.  We first make the change of variables $y=e^{\gamma c}Z_\epsilon(D_{h_\epsilon})$ to get
 \begin{align}\label{eqBbar2}
{B}_2(N,\epsilon)\leq &   C \,\E\Big[\mathbf{1}_{B_{N,\epsilon}}  \mathbf{1}_{I_{\epsilon}^c} 
  Z_\epsilon(D_{1}) ^{-\frac{\sigma }{\gamma}} \Big]
  \end{align}
We claim

\begin{lemma}\label{CDBM}
Let $B$ be a standard Brownian motion and $\beta>x>0$ and $\theta\in]0,1/2[$. Then, for some constant $C>0$ (independent of everything)
 \begin{align*}
 P_\epsilon(x):=&\P_x\Big(\min_{u\in[s_\epsilon,-\ln \epsilon]}\beta-B_u< s_\epsilon^\theta,\min_{u\in[0,-\ln \epsilon]}\beta-B_u\geq 0\Big)\\\leq &(\beta-x) (\ln\frac{_1}{^\epsilon})^{-1/2}s_\epsilon^{\theta-1/2 }.
 \end{align*}
\end{lemma}

Conditioning \eqref{eqBbar2} on the sigma algebra generated by $\{X_{\hat{g},u}(0);u>1\}$, we can use Lemma \ref{CDBM} to get
\begin{align}\label{eqBbar3}
(\ln\frac{_1}{^\epsilon})^\hf{B}_2(N,\epsilon)\leq &    Cs_\epsilon^{\theta-1/2 } \,\E\Big[( N-X_{\hat{g},1}(0))_+
  Z_\epsilon(D_{1}) ^{-\frac{\sigma }{\gamma}} \Big]
  \end{align}
 The last expectation is clearly finite and bounded independently of $\epsilon$ so that 
 \begin{equation}\label{B3}
\lim_{\epsilon\to 0}(\ln\frac{_1}{^\epsilon})^\hf{B}_2(N,\epsilon)=0.
\end{equation}
and, gathering \eqref{sumB}+\eqref{B1}+\eqref {B3limit}+\eqref{B3},
 the proof of  \eqref{est2} and hence Lemma \ref{crucial1} is complete.\qed

  \bigskip
  
\noindent {\it Proof of Lemma \ref{CDBM}.} We condition first on the filtration $\mathcal{F}_{s_\epsilon}$ generated by the Brownian motion up to time ${s_\epsilon}$. From Lemma \ref{standard}, we obtain
 \begin{align*}
 \P\Big[&\min_{u\in[s_\epsilon,-\ln \epsilon]}\beta-B_u< (s_\epsilon)^\theta,\min_{u\in[s_\epsilon,-\ln \epsilon]}\beta-B_u\geq 0|\mathcal{F}_{s_\epsilon}\Big]
 \\
 \leq &\sqrt{\frac{2}{\pi}} \int_{\frac{(\beta-B_{s_\epsilon}-(s_\epsilon)^\theta)_+}{(\ln\frac{1}{\epsilon}-s_\epsilon)^{1/2}}}^{\frac{(\beta-B_{s_\epsilon})_+}{(\ln\frac{1}{\epsilon}-s_\epsilon)^{1/2}}}e^{-\frac{u^2}{2}}\,du
 \\
  \leq &\mathbf{1}_{\{\beta-B_{s_\epsilon}\in [0,(s_\epsilon)^\theta]\}}\sqrt{\frac{2}{\pi}} \frac{(\beta-B_{s_\epsilon})_+}{(\ln\frac{1}{\epsilon}-s_\epsilon)^{1/2}} +\mathbf{1}_{\{\beta-B_{s_\epsilon}>(s_\epsilon)^\theta\}} \frac{(s_\epsilon)^\theta }{(\ln\frac{1}{\epsilon}-s_\epsilon)^{1/2}}.
   \end{align*}
Integrating, we get that
 \begin{align*}
 P_\epsilon(x)  \leq &\sqrt{\frac{2}{\pi}} (\ln\frac{1}{\epsilon}-s_\epsilon)^{-1/2} \E\Big[\mathbf{1}_{\{\min_{u\in[0,s_\epsilon]}\beta-B_u\geq 0\}}(\beta-B_{s_\epsilon})\mathbf{1}_{\{\beta-B_{s_\epsilon}\in [0,(s_\epsilon)^\theta]\}} \Big]\\
 &+ \frac{(s_\epsilon)^\theta }{(\ln\frac{1}{\epsilon}-s_\epsilon)^{1/2}}   \E\Big[\mathbf{1}_{\{\min_{u\in[0,s_\epsilon]}\beta-B_u\geq 0\}}\Big] .
   \end{align*}
 The second expectation is estimated with Lemma \ref{standard}. Concerning the first one, we use the fact under the probability measure $\frac{1}{\beta-x}\mathbf{1}_{\{\min_{u\in[0,s_\epsilon]}\beta-B_u\geq 0\}}(\beta-B_{s_\epsilon})$, the process $(\beta-B_u)_{u\leq s_\epsilon}$ is a $3d$-Bessel process, call it ${\rm Bess}_t$. Hence, using the Markov inequality, the scale invariance of a Bessel process and the fact that the mapping $x\mapsto \E^x\Big[\frac{1}{  {\rm Bess}_{1}}\Big] $ is decreasing, we deduce
  \begin{align*}
 P_\epsilon(x)  \leq &\sqrt{\frac{2}{\pi}} (\ln\frac{1}{\epsilon}-s_\epsilon)^{-1/2} (\beta-x)\E^{\beta-x}\Big[ \mathbf{1}_{\{{\rm Bess}_{s_\epsilon}\leq (s_\epsilon)^\theta\}} \Big] +\sqrt{\frac{2}{\pi}} \frac{(s_\epsilon)^{\theta-1/2 }}{(\ln\frac{1}{\epsilon}-s_\epsilon)^{1/2}} (\beta-x)  \\
  \leq &\sqrt{\frac{2}{\pi}} (\ln\frac{1}{\epsilon}-s_\epsilon)^{-1/2} (\beta-x)(s_\epsilon)^{\theta-1/2 }\E^0\Big[ \frac{1}{ {\rm Bess}_{1}}\Big] +\sqrt{\frac{2}{\pi}} \frac{(s_\epsilon)^{\theta-1/2 }}{(\ln\frac{1}{\epsilon}-s_\epsilon)^{1/2}} (\beta-x).\qed 
   \end{align*}

\appendix

\section{Auxiliary lemma}
 
 \begin{lemma}\label{standard}
We have for $\beta>0$
$$\P(\sup_{u\leq t}B_u\leq \beta)=\sqrt{\frac{2}{\pi}}\int_0^{\frac{\beta}{\sqrt{t}}} e^{- \frac{u^2}{2}}\,du\leq \sqrt{\frac{2}{\pi}} \frac{\beta}{\sqrt{t}}.$$
\end{lemma}  

The proof is elementary and thus left to the reader.

%
%

%

\end{document}